\documentclass[12pt]{amsart}

\usepackage{latexsym}
\usepackage{amssymb}
\usepackage{graphicx}
\usepackage{latexsym}
\usepackage{amssymb}

\newcommand{\excise}[1]{}

\newtheorem{thm}{Theorem}
\newtheorem{lemma}[thm]{Lemma}

\newtheorem{prop}[thm]{Proposition}

\newtheorem{Example}[thm]{Example}
\newtheorem{Remark}[thm]{Remark}
\newtheorem{Alg}[thm]{Algorithm}
\newtheorem{Defn}[thm]{Definition}

\newenvironment{remark}{\begin{Remark}\rm}
                {\mbox{}~\hfill$\square$\end{Remark}}

\newenvironment{defn}{\begin{Defn}\rm}
        {
         \end{Defn}}

    {\begin{eqnarray}}
    {\end{eqnarray}$\!\!$}

\newenvironment{eq*}%
    {\begin{eqnarray*}}
    {\end{eqnarray*}$\!\!$}

    {\begin{equation}}
    {\end{equation}$\!\!$}

\newenvironment{eqn*}%
    {\begin{equation*}}
    {\end{equation*}$\!\!$}

        {\begin{list}
                {\noindent\makebox[0mm][r]{\arabic{enumi}.}}
                {\leftmargin=5.5ex \usecounter{enumi}}
        }
        {\end{list}}

        {\begin{list}
                {\noindent\makebox[0mm][r]{(\roman{enumi})}}
                {\leftmargin=5.5ex \usecounter{enumi}}
        }
        {\end{list}}

\def\footrc#1{\hbox{\footnotesize${#1}$}}
\def\tinyrc#1{\hbox{\tiny${#1}$}}

\def\<{\langle}
\def\>{\rangle}
\def\0{{\mathbf 0}}
\def\1{{\mathbf 1}}
\def\BB{{\mathcal B}}
\def\CC{{\mathbb C}}

\def\FF{{\mathcal F}}
\def\GG{{\mathbf G}}

\def\LL{{\mathcal L}}

\def\OO{{\mathcal O}}
\def\PP{{\mathbb P}}

\def\TT{{\mathbf T}}
\def\ZZ{{\mathbb Z}}
\def\aa{{\mathbf a}}

\def\pp{{\mathbf p}}
\def\xx{{\mathbf x}}

\def\zz{{\mathbf z}}

\def\th{{\rm th}}
\def\conv{{\rm conv}}
\def\proj{{\rm Proj}}
\def\spec{{\rm Spec}}

\def\BD{B_\Delta}
\def\FL{{\mathcal F}{\ell}_n}
\def\IN{\mathsf{in}}
\def\UL{\Upsilon_{\!\lambda}}
\def\gl{{G_{\!}L}}
\def\mn{{M_n}}
\def\pg{P\dom\hspace{-.2ex}G}
\def\GIT{\dom\hspace{-.75ex}\dom\hspace{-.3ex}}
\def\dom{\backslash}
\def\gln{{G_{\!}L_n}}
\def\gli{{G_{\!}L_i}}
\def\glnone{{G_{\!}L_{n-1}}}
\def\too{\longrightarrow}
\def\bgln{B\dom\hspace{-.2ex}\gln}

\def\glnn{({G_{\!}L_n})^n}

\def\onto{\twoheadrightarrow}

\def\cprime{$'$}

\def\Sym{{\rm {Sym}}}

\def\ol#1{{\overline {#1}}}

\begin{document}

\title[Toric degeneration of Schubert varieties]%
      {Toric degeneration of Schubert varieties\\
    and Gel$'$fand--Cetlin polytopes}
\author{Mikhail Kogan}
\thanks{Both authors were supported by National Science Foundation
      Postdoctoral Research Fellowships}
\address{Northeastern University\\Boston, MA\\USA}
\email{misha@research.neu.edu}
\author{Ezra Miller}
\address{Mathematical Sciences Research Institute\\Berkeley, CA\\USA}
\email{emiller@msri.org}
\date{16 March 2003}

\begin{abstract}
\noindent
This note constructs the flat toric degeneration of the manifold $\FL$
of flags in~$\CC^n$ from \cite{GL96} as an explicit GIT quotient of the
Gr\"obner degeneration in~\cite{grobGeom}.  This implies that Schubert
varieties degenerate to reduced unions of toric varieties, associated to
faces indexed by rc-graphs (reduced pipe dreams) in the
Gel$'$fand--Cetlin polytope.  Our explicit description of the toric
degeneration of~$\FL$ provides a simple explanation of how
Gel$'$fand--Cetlin decompositions for irreducible polynomial
representations of~$\gln$ arise via geometric quantization.
\end{abstract}

\maketitle

{}


\section*{Introduction}
\label{sec:intro}

A number of recent developments at the intersection of algebraic
geometry and combinatorics have exploited degenerations of certain
varieties related to linear algebraic groups.  Sometimes the varieties
involved have been classical flag and Schubert varieties
\cite{KoganThesis,VakLR,VakSchub}, and other times they have been
closely related affine varieties \cite{grobGeom,KMS}.  In all cases, it
has been vital not only to construct an appropriate family degenerating
the primal variety, but also to identify combinatorially all components
occurring in the degenerate limit.  Indeed, this is what geometrically
produces (or reproduces) various combinatorial formulae (for classical
objects such as Littlewood--Richardson coefficients \cite{Ful97}, or for
universally defined polynomials discovered more recently
\cite{LS82a,LS82b,Ful99,BF,Buch02}): components of the special fiber
correspond to combinatorial summands in the desired formula.

Independently motivated degenerations of similar flavors have appeared
in areas related to representation theory, particularly standard
monomial theory \cite{GL96,Chi00,Cal02}.  The primal varieties there
have been generalized flag and Schubert varieties in arbitrary type,
with the limiting fibers being toric varieties, \mbox{or reduced unions
thereof}.

Our goal in this note is to create a single geometric framework relating
some of the more natural degenerations above.
To this end, we express the flat toric degeneration of the manifold
$\FL$ of flags in~$\CC^n$ from \cite{GL96}, which is a special case of
the degenerations in~\cite{Cal02}, as a quotient of the Gr\"obner
degeneration of $n \times n$ matrices~$\mn$ in~\cite{grobGeom}.  The
quotient is constructed by deforming the action of the lower triangular
matrices $B\subset \gln$ on the space~$\mn$ of matrices.  More
precisely, we construct an explicit action of $B$ on $\mn\times \CC$
and define the GIT quotient $B\GIT(\mn\times \CC)$, which is the total
family~$\FF$ of the desired degeneration and still fibers over the
line~$\CC$.

This degeneration can also be thought of simply as a Pl\"ucker embedding
of the Gr\"obner degeneration from~\cite{grobGeom}.  The closure of the
image of this embedding is the family~$\FF$ of projective varieties over
the affine line~$\CC$, whose fiber $\FF(1)$ over $1 \in \CC$ is the flag
variety~$\FL$, and whose fiber $\FF(0)$ over $0 \in \CC$ is the toric
variety associated to the {\em Gel\/$'\!$fand--Cetlin polytope}\/ from
representation theory \cite{GC50,GS83}.

Two consequences result from our explicit description of the
family~$\FF$.  First, since~$\FF$ is derived from Gr\"obner
degeneration, it induces a subfamily degenerating every Schubert variety
in~$\FL$.  Therefore---and this is the main point---applying the
combinatorial characterization of degenerated matrix Schubert varieties
in \cite{grobGeom}, we characterize in Theorem~\ref{t:main} {\em which
faces}\/ of the Gel$'$fand--Cetlin toric variety occur in degenerate
Schubert varieties of~$\FL$.  Namely, these faces correspond by
\cite{KoganThesis} to combinatorial diagrams called {\em rc-graphs}\/
(or {\em reduced pipe dreams}\/) \cite{FKyangBax,BB}.

Our second consequence is a simple explanation
(Section~\ref{sec:decomp}) for how the classical {\em
Gel\/$'\!$fand--Cetlin decomposition}\/ of irreducible polynomial
representations of~$\gln$ into one-dimensional weight spaces arises
geometrically from toric degeneration~$\FF$ of the flag manifold.
The idea is to think of Gel$'$fand--Cetlin decomposition as geometric
quantization on the total space of the family~$\FF$, directly extending
the manner in which the Borel--Weil theorem is geometric quantization at
the fiber $\FL = \FF(1)$.
The point is that sections of line bundles over~$\FL$, which constitute
irreducible representations of~$\gln$ by Borel--Weil, canonically
acquire at the special fiber of~$\FF$ an action of the torus densely
embedded in the toric variety~$\FF(0)$.
This produces a basis for sections over $\FL$ indexed by integer points
of the Gel$'$fand--Cetlin polytope.

\medskip
The methods in this note can be extended to partial flag manifolds in
type~$A$, but we believe the most exciting prospects for future research
lie in extensions to other types.  In particular, \cite{GL96} and
\cite{Cal02} describe a number of degenerations of generalized flag
varieties to toric varieties.  Under these degenerations, Schubert
varieties become unions of toric subvarieties.  In ``nice'' cases
\cite{Lit98} (this is a technical term), the degenerate toric variety
has an easily described moment polytope, in terms of generalizations of
Gel$'$fand--Cetlin patterns.  Identification of the components in
degenerations of Schubert varieties
could therefore provide arbitrary-type combinatorial generalizations of
pipe dreams,
such as those suggested by \cite{FKBn} in type~$B$.

The approach of degenerating generalized flag manifolds $\pg$
instead of degenerating groups~$G$ sidesteps the use of an
equivariant partial compactification of~$G$ to a vector space,
which is suggested by~\cite{grobGeom},
but is something we do not know how to do for arbitrary linear algebraic
groups.  Furthermore, one should be able to obtain positive
combinatorial formulae for Schubert classes in arbitrary type (though
perhaps not a definition of Schubert {\em polynomial}\/) by summing over
components in torically degenerated Schubert varieties.  The shapes
taken by such combinatorial formulae would echo the manner in which
\cite{KoganThesis} geometrically decomposes Schubert classes, agreeing
with the combinatorially positive formula in~\cite{BJS,FSnilCoxeter} for
the Schubert polynomials of Lascoux and Sch\"utzenberger \cite{LS82a}.

\medskip
Part of our purpose in writing this note was to make toric degenerations
of flag and Schubert varieties as in \cite{GL96,Chi00,Cal02} accessible
to an audience unaccustomed to specialized language surrounding
arbitrary finite type root systems and standard monomial theory.
In particular, we thought it important to include an equivalent
characterization of the family~$\FF$ in the language of sagbi bases
(Theorem~\ref{t:degeneration}), whose elementary definition and
properties we review in Section~\ref{sec:plucker}.  This version of the
toric degeneration of~$\FL$ is implicit in~\cite{GL96}, but so much so
as to be difficult to locate.  In addition, although the identification
of the limiting fiber~$\FF(0)$ as the Gel$'$fand--Cetlin toric variety
can be derived using results of \cite{Cal02} and~\cite{Lit98}, we
include the appropriate combinatorial arguments for
the reader's convenience (Section~\ref{sec:GC}).

\subsection*{Organization}
Our deformation of the Borel action on $\mn$ is constructed in
Section~\ref{sec:basics}.  Then Section~\ref{sec:plucker} uses the
results of \cite{GL96} to show that the quotient~$\FF =
B\GIT(\mn\times\CC)$ is a flat family.  Section~\ref{sec:GC} proves that
the zero fiber of~$\FF$ is the Gel$'$fand--Cetlin toric variety.  The
connection to~\cite{grobGeom} is exhibited in Section~\ref{sec:schubert}
by explaining how Schubert varieties behave inside of~$\FF$.  The final
section deals with geometric quantization of the degeneration, by
analyzing sections of line bundles over the family~$\FF$, thereby
constructing Gel$'$fand--Cetlin decompositions geometrically.

\section{Degenerating the Borel action}

\label{sec:basics}

Thinking about $\gln$ as a subset of $n\times n$ matrices~$\mn$ allows
us to think about the flag manifold $\FL = \bgln$ as a GIT quotient of
$\mn$ by $B$ whose precise definition is given in
Section~\ref{sec:plucker}.  In this section we construct a degeneration
of the action of $B$ on $\mn$, and in the next section we explain what
happens to the GIT quotient under this degeneration.

The group $\glnn$ has a left action on $\mn$ columnwise: if $Z \in \mn$
has columns $Z_1, \ldots, Z_n$, then $\gamma =
(\gamma_1,\ldots,\gamma_n) \in \glnn$ acts via
\begin{eqnarray} \label{gamma}
  \gamma Z &=& \hbox{matrix with columns } \gamma_1 Z_1, \ldots,
  \gamma_n Z_n.
\end{eqnarray}
The torus of $\glnn$ under the left action coincides with the standard
torus inside $\gl(\mn)$, scaling separately each entry of any given
matrix.

Let $\BD \subset \glnn$ be the image of the {\em lower}\/ triangular
Borel subgroup $B \subset \gln$ under the $n$-fold diagonal embedding in
$\glnn$, so $\BD = \{(b,\ldots,b) \mid b \in B\}$.

For every one-parameter subgroup $T \cong \CC^*$ inside the torus of
$\glnn$ consisting of sequences of diagonal matrices, denote by
$\tilde\tau \in T$ the element corresponding to the complex number~$\tau
\in \CC^*$.  Given a matrix $\omega = (\omega_{ij})$ of integers, let
the one-parameter subgroup $T(\omega)$ consist of sequences of diagonal
matrices, with the $j^\th$ component of~$\tilde\tau$ being the diagonal
matrix $\tilde\tau_{\!j} = {\rm{diag}} (\tau^{\omega_{1j}},\ldots,
\tau^{\omega_{nj}})$ for the $j^\th$ column~of~$\omega$.

For the rest of the paper, fix the
matrix~$\omega$ whose entries equal
%
\begin{equation} \label{eq:matrix}
\begin{array}{r@{\ }c@{\ }l@{\ }l}
                \omega_{ij}&=&3^{n-i-j}&\text{ if } i+j\leq n, \\
\hbox{and\quad} \omega_{ij}&=&0&\text{ if }i+j>n.
\end{array}
\end{equation}
For instance, when $n=5$,
we get the following $5\times 5$ matrix:
\begin{eqnarray*}
\omega &=&
\footrc{
\left[\!\!\begin{array}{ccccc}
 27 & 9 & 3 & 1 & 0 \\
 9 & 3 & 1 & 0 & 0 \\
 3 & 1 & 0 & 0 & 0 \\
 1 & 0 & 0 & 0 & 0 \\
 0 & 0 & 0 & 0 & 0
\end{array}\!\right]}
\end{eqnarray*}

Consider the family $\BB^* \subset B^n \times \CC^*$ of subgroups
of~$B^n$ with fiber
\begin{eqnarray} \label{tau}
  B(\tau) &=& \tilde\tau^{-1} \BD \tilde\tau
\end{eqnarray}
over $\tau \in \CC^*$, where $\tilde\tau$ lies in the one-parameter
subgroup $T(\omega)$ corresponding to~$\omega$.  When $n=5$, for
instance, this multiplies the entries in $B^n$ by powers of~$\tau$ as
follows.
\begin{eqnarray*}
&
\def\*#1{\makebox[.4ex][c]{$\tau^{#1}$}}
\def\lt{\multicolumn{1}{|c}{\makebox[.4ex][c]{$\tau$}}}
\def\mt{\multicolumn{1}{c}{\makebox[.4ex][c]{$\tau$}}}
\def\mc{\multicolumn{1}{|c}{}}
\def\mo{\multicolumn{1}{c}{\makebox[.4ex][c]{$1$}}}
\def\lo{\multicolumn{1}{|c}{\makebox[.4ex][c]{$1$}}}
\def\ro{\multicolumn{1}{c|}{\makebox[.4ex][c]{$1$}}}
\def\mss{\multicolumn{1}{|c|}{\makebox[.4ex][c]{$1$}}}
\def\mcc{\multicolumn{1}{|c|}{}}
\left(\ \tinyrc{\begin{array}{|ccccc|}\hline
    \lo     &\mc     &        &        &    \\\cline{1-2}
   \ \*{18}     &\lo     &\mc     &        &    \\\cline{2-3}
   \  \*{24}     &\*6     &\lo     &\mc     &    \\\cline{3-4}
   \ \*{26}     &\*8     &\*2     &\lo     &\mcc\\\cline{4-5}
   \ \*{27}     &\*9     &\*3     &\mt     &\mss\\\hline
\end{array}}
\ ,\
\tinyrc{\begin{array}{|ccccc|}\hline
    \lo     &\mc     &        &        &    \\\cline{1-2}
    \*6     &\lo     &\mc     &        &    \\\cline{2-3}
    \*8     &\*2     &\lo     &\mc     &    \\\cline{3-4}
    \*9     &\*3     &\mt     &\lo     &\mcc\\\cline{5-5}
    \*9     &\*3     &\mt     &\lo     &\ro \\\hline
\end{array}}
\ ,\
\tinyrc{\begin{array}{|ccccc|}\hline
    \lo     &\mc     &        &        &    \\\cline{1-2}
    \*2     &\lo     &\mc     &        &    \\\cline{2-3}
    \*3     &\mt     &\lo     &\mc     &    \\\cline{4-4}
    \*3     &\mt     &\lo     &\mo     &\mcc\\\cline{5-5}
    \*3     &\mt     &\lo     &\mo     &\ro \\\hline
\end{array}}
\ ,\
\tinyrc{\begin{array}{|ccccc|}\hline
    \lo     &\mc     &        &        &    \\\cline{1-2}
    \lt     &\lo     &\mc     &        &    \\\cline{3-3}
    \lt     &\lo     &\mo     &\mc     &    \\\cline{4-4}
    \lt     &\lo     &\mo     &\mo     &\mcc\\\cline{5-5}
    \lt     &\lo     &\mo     &\mo     &\ro \\\hline
\end{array}}
\ ,\
\tinyrc{\begin{array}{|ccccc|}\hline
    \lo    &\mc     &        &        &    \\\cline{2-2}
    \lo    &\mo     &\mc     &        &    \\\cline{3-3}
    \lo    &\mo     &\mo     &\mc     &    \\\cline{4-4}
    \lo    &\mo     &\mo     &\mo     &\mcc\\\cline{5-5}
    \lo    &\mo     &\mo     &\mo     &\ro \\\hline
\end{array}}
\ \right)
\end{eqnarray*}
The family $\BB^*$ extends to a family over all of~$\CC$:

\begin{defn}
The family $\BB \subset B^n \times \CC$ has fiber $B(\tau)$ over $\tau
\in \CC^*$, and fiber $B(0)$ consisting of sequences $(b_1,\ldots,b_n)
\in B^n$, where $b_j$ is obtained from the matrix $b_n$ by setting
to~$0$ all entries in columns $1,\ldots, n-j$ strictly below the main
diagonal.
\end{defn}

When $n=5$, elements in the special fiber $B(0)$ look heuristically
like:
\begin{eqnarray*}
&
\def\*{\makebox[.4ex][c]{$b$}}
\def\mc{\multicolumn{1}{|c}{}}
\def\ms{\multicolumn{1}{|c}{\*}}
\def\mss{\multicolumn{1}{|c|}{\*}}
\def\mcc{\multicolumn{1}{|c|}{}}
\left(\
\tinyrc{\begin{array}{|ccccc|}\hline
    \ms&\mc&   &   &    \\\cline{1-2}
       &\ms&\mc&   &    \\\cline{2-3}
       &   &\ms&\mc&    \\\cline{3-4}
       &   &   &\ms&\mcc\\\cline{4-5}
       &   &   &   &\mss\\\hline
\end{array}}
\ ,\
\tinyrc{\begin{array}{|ccccc|}\hline
    \ms&\mc&   &   &    \\\cline{1-2}
       &\ms&\mc&   &    \\\cline{2-3}
       &   &\ms&\mc&    \\\cline{3-4}
       &   &   &\ms&\mcc\\\cline{5-5}
       &   &   &\ms&\*  \\\hline
\end{array}}
\ ,\
\tinyrc{\begin{array}{|ccccc|}\hline
    \ms&\mc&   &   &    \\\cline{1-2}
       &\ms&\mc&   &    \\\cline{2-3}
       &   &\ms&\mc&    \\\cline{4-4}
       &   &\ms&\* &\mcc\\\cline{5-5}
       &   &\ms&\* &\*  \\\hline
\end{array}}
\ ,\
\tinyrc{\begin{array}{|ccccc|}\hline
    \ms&\mc&   &   &    \\\cline{1-2}
       &\ms&\mc&   &    \\\cline{3-3}
       &\ms&\* &\mc&    \\\cline{4-4}
       &\ms&\* &\* &\mcc\\\cline{5-5}
       &\ms&\* &\* &\*  \\\hline
\end{array}}
\ ,\
\tinyrc{\begin{array}{|ccccc|}\hline
    \*&\mc&   &   &    \\\cline{2-2}
    \*&\* &\mc&   &    \\\cline{3-3}
    \*&\* &\* &\mc&    \\\cline{4-4}
    \*&\* &\* &\* &\mcc\\\cline{5-5}
    \*&\* &\* &\* &\*  \\\hline
\end{array}}
\ \right)
\end{eqnarray*}

\begin{lemma} \label{l:tau}
There is a canonical algebraic group isomorphism $B \times \CC \to
\BB$ over~$\CC$.
\end{lemma}
\begin{proof}
Use that $B \cong \BD$ by sending $b \mapsto b_\Delta = (b,\ldots,b)$.
For $\tau \neq 0$ the isomorphism is now by~(\ref{tau}), sending
$b_\Delta$ to $\tilde\tau^{-1} b_\Delta \tilde\tau$.  For $\tau = 0$,
the map sets to~$0$ all entries in columns $1,\ldots, n-j$ strictly below
the main diagonal in the $j^\th$ entry of $(b,\ldots,b)$.%

Elementary computation shows that if~$b$ is a lower triangular matrix,
then the matrix $\tilde\tau_j^{-1} b \tilde\tau_j$ has no negative
powers of~$\tau$, and setting $\tau$ to~$0$ has the effect of setting
to~$0$ all entries in columns $1,\ldots, n-j$ strictly below the main
diagonal of~$b$.  Hence the $\tau = 0$ case above really is obtained
from the $\tau \neq 0$ case by taking limits as $\tau \to 0$.%
\end{proof}

The family $\BB$ of groups acts fiberwise on~$\mn \times \CC$, but
Lemma~\ref{l:tau} allows us to view this fiberwise action as a single
action of~$B$ on the total space $\mn \times \CC$.  The actions on all
fibers $\mn \times \tau$ are isomorphic for $\tau \in \CC^*$, in the
sense that the map $Z \times 1 \mapsto \tilde\tau^{-1} Z \times \tau$
identifies $\mn \times 1$ with $\mn \times \tau$ equivariantly with
respect to the actions of~$B$ on the fibers over~$1$ and~$\tau$.
However, when $\tau$ equals zero, $B$~acts on the $j^\th$ column as the
product of an $n-j$ dimensional torus (in the upper-left corner) and a
smaller Borel group with $j$~columns
(in the lower-right corner).

The action of $B = B(0)$ on $\mn\times 0$ commutes with an $\binom n2$
dimensional torus action, which scales all entries lying strictly above
the main antidiagonal in each $n\times n$ matrix.  We shall see that
this torus acts on the degenerated $\FL$ to make it a toric variety.

\section{Degeneration of Pl\"ucker coordinates}

\label{sec:plucker}

Degenerating the action of $B$ on $\mn$ via the action of $B \times
\CC$ in Lemma~\ref{l:tau} on $\mn \times \CC$ induces a degeneration of
the GIT quotient of~$\mn$ by~$B$.
Using the results of Gonciulea and Lakshmibai \cite{GL96}, we show
that the GIT quotient $B\GIT (\mn\times \CC)$, when defined
appropriately, flatly degenerates the flag manifold~$\FL$ to a toric
variety.

Let $U$ be the lower triangular matrices with $1$'s on the diagonal (the
unipotent radical of~$B$).
As can be done in the general setting of $B$ actions, we define the GIT
quotient of~$\mn$ by $B$ to be the ``multiple $\proj$'' of the ring of
$U$-invariant functions on~$\mn$.  Let us be more precise in the present
case.

For a subset $J\subset\{1,\ldots,n\}$ of size~$k$, define
$\Delta_J(Z)$ to be the minor of an $n\times n$ matrix~$Z$ whose
columns are given by the set~$J$ and whose rows $1,\ldots,k$ are
top-justified.  Writing $\CC[\zz]$ with $\zz = (z_{ij})_{i,j=1}^n$
for the coordinate ring of~$\mn$, the set of {\em Pl\"ucker
coordinates}\/ consists of all minors having the form
\begin{eqnarray*}
  p_J &=& \Delta_J(n \times n \hbox{ matrix of variables $\zz$}).
\end{eqnarray*}
They generate the ring $\CC[\pp] \subset \CC[\zz]$ of $U$ invariant
functions on~$\mn$.  This invariant ring~$\CC[\pp]$ can be expressed as
a quotient
\begin{eqnarray*}
  \CC[\xx^1] \otimes \cdots \otimes \CC[\xx^n] &\onto& \CC[\pp]
\end{eqnarray*}
of the tensor product over~$\CC$ of~$n$ polynomial rings
$\CC[\xx^k]$, where $\xx^k$ is a set of variables~$x_J$ indexed by
the size~$k$ subsets of $\{1,\ldots,n\}$.
Thus the spectrum of~$\CC[\pp]$ is a subvariety of $\bigwedge^{\!*}
\CC^n$.  This gives rise to the {\em multiple}\/ $\proj$ of~$\CC[\pp]$,
by which we mean the corresponding subscheme of~$\prod^n_{k=1}
\PP(\bigwedge^k\CC^n) = \prod^n_{k=1}\proj(\CC[\xx^k])$.  This {\em
Pl\"ucker embedding}\/ of the flag manifold~$\FL$ is the GIT quotient
of~$\mn$ by~$B$.

Now let us turn to the GIT quotient of $\mn \times \CC$ by the action
of~$B$ constructed in the previous section.  In this context, we are
thinking of $\mn \times \CC$ as a (trivial) family over~$\CC$, and we
wish to quotient out by the fiberwise action of the family $\BB$ of
groups parametrized by~$\CC$.  By definition, this GIT quotient is the
multiple $\proj$ of the $\CC[t]$-algebra of
$U$-invariant functions on $\mn \times \CC$.
To describe some of these invariant functions, we need a preliminary
result.

Think of the matrix $\omega$ as a weighting on the coordinate ring
$\CC[\zz]$ of~$\mn$ under which each variable~$z_{ij}$ has
weight~$\omega_{ij}$.  Also, let $\Delta_{I,J}$ take the minor with
rows~$I$ and columns~$J$.  The {\em antidiagonal}\/ of such a minor is
the product of all entries along the main antidiagonal in the
corresponding square matrix.

\begin{lemma} \label{l:weights}
If every variable dividing the antidiagonal term of the minor
$\Delta_{I,J}(Z)$ in the generic matrix~$Z$ lies on or above the main
antidiagonal of~$Z$, then the unique \mbox{$\zz$-monomial} in
$\Delta_{I,J}(Z)$ with the lowest weight is its antidiagonal term.
\end{lemma}
\begin{proof}
It suffices to prove the lemma when $I$ and $J$ have cardinality~$2$,
because every $\zz$-monomial in each minor can be made into the
antidiagonal term by successively replacing $2 \times 2$ diagonals with
$2 \times 2$ antidiagonals. In the $2 \times 2$ case, let $I = \{i,
i+k\}$ and $J = \{j, j+\ell\}$ with $i,j,k,\ell \geq 1$.  The weights on
the two terms in $\Delta_{I,J}(Z)$ satisfy
\begin{eq*}
  3^{n-i-j} +3^{n-i-k-j-\ell} > 3^{n-i-k-j} +3^{n-i-j-\ell},
\end{eq*}
which proves the lemma.
\end{proof}

Denote by~$\tilde t$ the $n$-tuple of $n\times n$ diagonal matrices
whose $j^\th$ diagonal entry in the $i^\th$ matrix is~$t^{\omega_{ij}}$,
and define $\tilde t Z$ for $\tilde t = \gamma$ as in~(\ref{gamma}) for
the matrix~$Z$ of variables.
In addition, for $J=\{j_1,\cdots,j_k\}$, let
\begin{eqnarray*}
  \omega_J &=& \sum_{i=1}^k\omega_{i,n+1-j_i}
\end{eqnarray*}
be the sum of weights along the antidiagonal of the square submatrix in
rows $1,\ldots,k$ and columns~$J$ of~$\omega$.  Then, as an immediate
consequence of Lemma~\ref{l:weights}, we conclude that the polynomials
\begin{eqnarray} \label{q}
  q_J &=& t^{-\omega_J}\Delta_J(\tilde t Z)
\end{eqnarray}
are $U$-invariants in $\CC[\zz,t]$ under the action of~$B$ resulting
from Lemma~\ref{l:tau}.  The power $t^{-\omega_J}$ precisely makes the
antidiagonal term of~$q_J$ have coefficient~$\pm 1$.

\begin{defn} \label{d:FF}
Define the family~$\FF$ inside the product $\prod^n_{k=1}
\PP(\bigwedge^k\CC^n)\times \CC$ over the line $\CC = \spec(\CC[t])$ as
the multiple $\proj$ of the subalgebra \mbox{$\CC[q_J \mid J \subseteq
\{1,\ldots,n\}]$ of~$\CC[\zz,t]$}.
\end{defn}

We wish to state the main result in this section in terms of sagbi
bases.  Recall that a {\em term order}\/ on $\CC[\zz]$ is a
multiplicative total order on monomials with $1 \in \CC[\zz]$ being
smaller than any other monomial; see \cite[Chapter~15]{Eis}.  A~set
$\{f_1,\ldots,f_r\} \subset \CC[\zz]$ is a {\em sagbi basis}\/ if the
initial term $\IN(f)$ of every polynomial~$f$ in the subalgebra
$\CC[f_1,\ldots,f_r]$ lies inside the {\em initial subalgebra}\/
generated by the initial terms $\IN(f_1),\ldots,\IN(f_r)$.  The initial
subalgebra is generated by monomials, so its multiple $\proj$ is a toric
variety.  Choosing a weight order inducing the given term order
\cite[Chapter~15]{Eis} allows us to express the original algebra and its
initial algebra as the fibers over~$0$ and~$1$ of a flat family of
subalgebras of~$\CC[\zz]$.  In fact, $\{f_1,\ldots,f_r\}$ form a sagbi
basis if and only if this degeneration to the initial subalgebra is
flat.

The terms orders on the coordinate ring~$\CC[\zz]$ of~$\mn$ that
interest us are {\em antidiagonal}\/ and {\em diagonal}.  By~definition,
the leading term of any minor in the matrix of variables under an
(anti)diagonal term order is its (anti)diagonal term, namely the product
of all entries on the (anti)diagonal of the corresponding square
submatrix.  Initial terms of polynomials other than minors will not be
important in what follows.

\begin{thm} \label{t:degeneration}
The polynomials~$q_J$ from~(\ref{q}) generate the $\CC[t]$-algebra of
$U$-invariant functions inside~$\CC[\zz,t]$, so $\FF$ is the GIT
quotient family $B\GIT (\mn \times \CC)$ flatly degenerating the flag
manifold $\FL=\FF(1)$ to a toric variety~$\FF(0)$.  In fact, the
Pl\"ucker coordinates~$p_J$ constitute a sagbi basis for any diagonal or
antidiagonal term~order.
\end{thm}
\begin{proof}
The second sentence implies the first;
to show this, we may as well assume by symmetry (reflecting
left-to-right) that the term order is antidiagonal.

Setting $t = 1$ in~(\ref{q}) obviously yields the Pl\"ucker
coordinates~$p_J$.  Therefore the polynomials~$q_J$ generate the
$U$-invariants in $\CC[\zz,t^{\pm1}]$ over the coordinate
ring~$\CC[t^{\pm1}]$ of~$\CC^*$, because the Pl\"ucker coordinates~$p_J$
generate the $U$-invariants over \mbox{$t = 1$}, and the family of
$U$-invariants is trivial by scaling outside $t = 0$.  Intersecting the
$U$-invariants in $\CC[\zz,t^{\pm1}]$ with $\CC[\zz,t]$ yields
$U$-invariants in~$\CC[\zz,t]$, because a polynomial function on $\mn
\times \CC$ is $U$-invariant if and only if its restriction to~$\mn
\times \CC^*$ is.
Therefore, we must show that the polynomials~$q_J$ generate as a
$\CC[t]$ algebra the intersection with~$\CC[\zz,t]$ of the subalgebra
they generate inside $\CC[\zz,t^{\pm1}]$.  This follows from the sagbi
property.

For the proof of the second sentence, we assume the term order is
diagonal, to agree with \cite{GL96}.  Let $H$ be the set of all subsets
of $\{1,\ldots,n\}$.  Following \cite{GL96}, define a partial order on
$H$ as follows.  For $I=\{i_1<\cdots <i_k\}$ and $J=\{j_1<\cdots<
j_\ell\}$ set
\begin{eqnarray*}
  I\geq J &\Longleftrightarrow& k\leq \ell\text{ and } i_s\geq j_s
  \hbox{ for } 1 \leq s \leq k.
\end{eqnarray*}
This makes $H$ a distributive lattice.  It is shown in \cite{GL96} that
if $k\leq \ell$, then the meet and join of~$I$ and~$J$ are characterized
by
\begin{eqnarray*}
  I \wedge J &=& (\min(i_1,j_1),\ldots, \min(i_{k},j_{k}),
  j_{\ell+1},\ldots, j_{\ell})
\\
\hbox{and\ \ }
  I\vee J &=& (\max(i_1,j_1),\ldots, \max(i_{k},j_{k})).
\end{eqnarray*}
In the degeneration
from \cite[Theorems~5.2 and~10.6]{GL96}, the fiber over~$1$ equals the
subalgebra of~$\CC[\zz]$ generated by Pl\"ucker coordinates.  The
initial algebra in \cite{GL96} maps surjectively onto the `diagonal'
semigroup algebra generated by he diagonals, because the degenerated
Pl\"ucker relations $x_Ix_J - x_{I\wedge J}x_{I \vee J}$ in \cite{GL96}
hold on diagonals of Pl\"ucker coordinates. But the initial algebra in
\cite{GL96} and the diagonal semigroup algebra are integral domains of
equal Krull dimension~$\binom {n+1}2$ so the surjection must be an
isomorphism.  It follows that the diagonal semigroup algebra has the
same Hilbert series as the Pl\"ucker algebra, by flatness of the
degeneration in~\cite{GL96}.  Since the diagonal semigroup algebra is
certainly contained inside the initial algebra, whose Hilbert series
equals that of the Pl\"ucker algebra by flatness again, the diagonal
semigroup algebra must equal the initial algebra.
\end{proof}

\begin{remark}
For each nonzero $\tau$, the fiber $\FF(\tau)$ equals
$B(\tau)\hspace{-.2ex}\dom\hspace{-.2ex}\gln(\tau)$, where $\gln(\tau)$
is the set of $n\times n$ matrices~$Z$ with nonzero $\tau$-determinant
$\det(\tilde\tau Z)$.  When $\tau=0$, we can define $\gln(0)$ as the set
of matrices with nonzero entries along the antidiagonal.  The quotient
$B(0)\hspace{-.2ex}\dom\hspace{-.2ex}\gln(0)$ is still well defined, but
no longer equal to~$\FF(0)$: geometrically, under the zeroth Pl\"ucker
map, whose coordinates are obtained by setting $t = 0$ in~(\ref{q}), the
variety $\gln(0)$ has image equal to an open cell inside the toric
variety~$\FF(0)$.
\end{remark}

\section{Gel\cprime fand--Cetlin polytopes}

\label{sec:GC}

Next we show that the zero fiber $\mathcal F(0)$ is the toric variety
associated to a Gel$'$fand--Cetlin polytope, defined as follows.  Let
$\lambda = (\lambda_1>\cdots>\lambda_n)$ be a nonincreasing sequence of
nonnegative integers.  An array $\Lambda = (\lambda_{i,j})_{i+j\leq
n+1}$ of real numbers is a {\em Gel\/$'\!$fand--Cetlin pattern}\/
for~$\lambda$ if $\lambda_{i,1} = \lambda_i$ for all $i = 1,\ldots,n$,
and $\lambda_{i,j} \geq \lambda_{i,j+1} \geq \lambda_{i+1,j}$ for $i,j =
1,\ldots,n$.  Equivalently, entries in Gel$'$fand--Cetlin
patterns~$\Lambda$ decrease in the directions indicated by the arrows in
diagram below, whose left column is~$\lambda$:
\begin{equation} \label{eq:GC}
\begin{array}{c@{\ }c@{\ }c@{\ }c@{\ }c@{\ }c@{\ }c}
\lambda_{1,1} & \to & \lambda_{1,2} & \to & \lambda_{1,3} & \to & \cdots
\\
\downarrow&\swarrow&\downarrow&\swarrow&
\\
\lambda_{2,1} & \to & \lambda_{2,2} & \to & \cdots
\\
\downarrow    &\swarrow&
\\
\lambda_{3,1} & \to & \cdots
\\
\downarrow
\\
\vdots
\end{array}
\end{equation}
The {\em Gel\/$'\!$fand--Cetlin polytope}\/ $P_\lambda$ is the convex
hull of all integer Gel$'$fand--Cetlin patterns for~$\lambda$.
This polytope defines the {\em Gel\/$'\!$fand--Cetlin toric variety}\/
together with its projective embedding.  For background on toric
varieties, see \cite{Ful93}.

Set $a_k=\lambda_k-\lambda_{k+1}$ for $k=1,\ldots, n$, where by
convention \mbox{$\lambda_{n+1}=0$}, and assume $a_k \geq 1$ for
all~$k$.  Recall that we expressed the flag manifold~$\FL$ as a
subvariety of the product \mbox{$\PP_1 \times \cdots \times \PP_n$},
where $\PP_k = \PP(\bigwedge^{\!k}\CC^n)$.  Since all~$a_k$ are strictly
positive, the integer sequence $\aa = (a_1,\ldots,a_n)$ corresponds to a
choice of very ample line bundle $\OO_{\FL}(\aa)$ on~$\FL$, namely the
result of tensoring together the pullbacks of the bundles
$\OO_{\PP_1}(a_1), \ldots, \OO_{\PP_n}(a_n)$ to~$\FL$.  In fact, we get
a choice of very ample line bundle $\OO_\FF(\aa)$ on the entire family
$\FF$ from Definition~\ref{d:FF}, to get an embedding of the family
\begin{eqnarray} \label{PL}
  \FF \to \PP_\lambda \times \CC &\hbox{where}& \PP_\lambda =
  \PP\big(\!\otimes^n_{k=1} \Sym^{a_k}({\textstyle
  \bigwedge^{\!k}}\CC^n)\big).
\end{eqnarray}
The image of the zero fiber $\FF(0)$ is a toric variety projectively
embedded inside $\PP_\lambda$ by a line bundle~$\OO_{\FF(0)}(\aa)$.

Let $\alpha_I$ be the exponent vector on the antidiagonal monomial
of the Pl\"ucker coordinate~$p_I$.  Such exponent vectors are elements
in~$\ZZ^{n^2}$ that look, for example, like
\begin{eqnarray*}
\def\*{\makebox[.4ex][c]{$1$}}
\tinyrc{\begin{array}{|@{\:\;}c@{\:\;}|@{\:\;}c@{\:\;}|@{\:\;}c@{\:\;}|%
                        @{\:\;}c@{\:\;}|}\hline
      &  &  &\*\\\hline
      &\*&  &  \\\hline
    \*&  &  &  \\\hline
      &  &\,&  \\\hline
\end{array}}\ \hbox{ for } p_{124} &\quad\hbox{and}\quad&
\def\*{\makebox[.4ex][c]{$1$}}
\tinyrc{\begin{array}{|@{\:\;}c@{\:\;}|@{\:\;}c@{\:\;}|@{\:\;}c@{\:\;}|%
                        @{\:\;}c@{\:\;}|}\hline
      &  &\*&  \\\hline
    \*&  &  &  \\\hline
      &\,&  &  \\\hline
      &  &  &\,\\\hline
\end{array}}\ \hbox{ for } p_{13}.
\end{eqnarray*}
The vector space of global sections of~$\OO_{\FF(0)}(\aa)$ decomposes
into a multiplicity-one direct sum of weight spaces.  The set of weights
occurring in this decomposition is
\begin{eqnarray} \label{UL}
  \quad \UL &=& \bigl\{\hbox{sums $\sum \alpha_I$ in which } a_k \hbox{
  of the indices } I \subseteq \{1,\ldots,n\} \hbox{ have size }
  k\bigr\}.
\end{eqnarray}

\begin{prop} \label{prop:GC-polytope}
The projective embedding $\FF(0)\to \PP_\lambda$ is the projective
embedding of the Gel\/$'\!$fand--Cetlin toric variety associated to the
polytope $P_\lambda$.
\end{prop}
\begin{proof}
By standard results about projective embeddings of toric varieties
\cite{Ful93}, it suffices to identify $P_\lambda$ with the convex hull
$\conv(\UL)$ of all points in~$\UL$.  This we shall prove for
all~$\lambda$, not just those producing strictly
positive~$a_1,\ldots,a_n$.

Denote the set of integer points of the polytope $P_\lambda$ by
$\Pi_\lambda$. Then $\Pi_\lambda$ sits inside an integer lattice
of rank $\frac{n(n-1)}{2}$, and $P_\lambda = \conv(\Pi_\lambda)$.

We claim that the linear map $\phi:\ZZ^{n^2} \to \ZZ^{\frac{n(n-1)}{2}}$
given by
\begin{eqnarray*}
  \lambda_{ij} &=& a_{i,j} + a_{i,j+1} + \cdots + a_{i,n+1-i}
\end{eqnarray*}
provides a bijection between $\UL$ and~$\Pi_\lambda$, implying the
required identification of $\conv(\UL)$ with~$P_\lambda$.

To check this, consider the map $\psi:\ZZ^{\frac{n(n-1)}{2}}\to
\ZZ^{n^2}$ given by
\begin{eqnarray*}
  a_{ij} &=& \lambda_{i,j}-\lambda_{i-1,j}.
\end{eqnarray*}
Notice that the composite maps $\phi\circ\psi$ and $\psi\circ\phi$ are
identities on $\Pi_\lambda$ and $\UL$ respectively.  It remains to check
$\phi(\UL)\subseteq \Pi_\lambda$ and $\psi(\Pi_\lambda)\subseteq \UL$;
these are simple exercises in linear algebra that go as follows.

To check $\phi(\UL)\subseteq \Pi_\lambda$, consider an element
$\alpha=\sum \alpha_I$ of $\UL$.  The $\lambda_{i,1}$ coordinate
of $\phi(\alpha)$ is the sum of all $a$'s in row~$i$, and this
equals $a_i+\cdots +a_n=\lambda_i$, which is the number of
indices~$I$ of size at least~$i$.  The $\lambda_{i,j}$ coordinate
of $\phi(\alpha)$ is the sum of the $a$'s in the horizontal strip
between $a_{i,j}$ and the antidiagonal. This is not greater
than~$\lambda_{i-1,j}$, which is the sum of the one-longer
horizontal strip starting at $a_{i-1,j}$. On the other hand,
$\lambda_{i,j}$ is at least $\lambda_{i-1,j+1}$. Indeed,
$\lambda_{i-1,j+1}$ is the sum of the entries in the horizontal
strip of the same length as for $\lambda_{i,j}$, but shifted down
one row and moved one column to the left. Since $\alpha$ is the
sum of $\alpha_I$'s, the sum of entries for $\lambda_{i-1,j+1}$ is
at most that for $\lambda_{i,j}$.  Thus $\phi(\UL)$ is a subset
of~$\Pi_\lambda$.

Conversely, given a Gel$'$fand--Cetlin pattern $\Lambda$ for
$\lambda$ we need to show that $\psi(\lambda)$ can be written as a
sum $\sum \alpha_I$.  (It will follow immediately that there are
$a_k$ indices~$I$ of cardinality $k$, since the number of indices
of cardinality at least $k$ must be $\lambda_k$.)  Derive a
pattern $\Lambda'$ from~$\Lambda$ by decreasing the last nonzero
entry $\lambda_{k,i_k}$ in each row by~$1$. Then $\Lambda'$ is
still a Gel$'$fand--Cetlin pattern.  At the same time it is clear
that $\psi(\Lambda)=\psi(\Lambda')+\alpha_I$ for the set
$I=\{i_1>\ldots>i_\ell\}$, where row~$\ell$ of~$\Lambda$ is not
zero but row~$\ell+1$ of~$\Lambda$ is zero.  Induction on the sum
of the entries
in Gel$'$fand--Cetlin patterns finishes the proof.%
\end{proof}

\section{Degenerating Schubert varieties}

\label{sec:schubert}

In this section we present our main theorem, which says that Schubert
varieties degenerate inside the family~$\FF$ to unions of toric
subvarieties given by rc-faces of the Gel$'$fand--Cetlin polytope.
First, we must review the combinatorics involved.

\begin{figure}[ht]
\small
\begin{picture}(300,90)
\put(0,64){1} \put(0,48){2} \put(0,32){3} \put(0,16){4}
\put(0,0){5}

\put(18,83){1} \put(34,83){2} \put(50,83){3} \put(66,83){4}
\put(82,83){5}

\put(81,71){\oval(8,8)[br]} \put(65,55){\oval(8,8)[br]}
\put(49,39){\oval(8,8)[br]} \put(33,23){\oval(8,8)[br]}
\put(17,7){\oval(8,8)[br]}

\put(17,71){\oval(8,8)[br]} \put(25,63){\oval(8,8)[tl]}

\put(33,71){\oval(8,8)[br]} \put(41,63){\oval(8,8)[tl]}

\put(49,71){\oval(8,8)[br]} \put(57,63){\oval(8,8)[tl]}

\put(65,71){\oval(8,8)[br]} \put(73,63){\oval(8,8)[tl]}

\put(17,23){\oval(8,8)[br]} \put(25,15){\oval(8,8)[tl]}



\put (18,35){\line(1,0){6}} \put (21,32){\line(0,1){6}}

\put (18,51){\line(1,0){6}} \put (21,48){\line(0,1){6}}


\put (34,35){\line(1,0){6}} \put (37,32){\line(0,1){6}}

\put (34,51){\line(1,0){6}} \put (37,48){\line(0,1){6}}


\put (50,51){\line(1,0){6}} \put (53,48){\line(0,1){6}}




\put (08,67){\line(1,0){10}} \put (08,51){\line(1,0){10}} \put
(08,35){\line(1,0){10}} \put (08,19){\line(1,0){10}} \put
(08,3){\line(1,0){10}} \put (24,67){\line(1,0){10}} \put
(24,51){\line(1,0){10}} \put (24,35){\line(1,0){10}} \put
(24,19){\line(1,0){10}} \put (40,67){\line(1,0){10}} \put
(40,51){\line(1,0){10}} \put (40,35){\line(1,0){10}} \put
(56,67){\line(1,0){10}} \put (56,51){\line(1,0){10}} \put
(72,67){\line(1,0){10}}


\put (21,70){\line(0,1){10}} \put (21,54){\line(0,1){10}} \put
(21,38){\line(0,1){10}} \put (21,22){\line(0,1){10}} \put
(21,6){\line(0,1){10}} \put (37,70){\line(0,1){10}} \put
(37,54){\line(0,1){10}} \put (37,38){\line(0,1){10}} \put
(37,22){\line(0,1){10}} \put (53,70){\line(0,1){10}} \put
(53,54){\line(0,1){10}} \put (53,38){\line(0,1){10}} \put
(69,70){\line(0,1){10}} \put (69,54){\line(0,1){10}} \put
(85,70){\line(0,1){10}}



\put(100,64){1} \put(100,48){2} \put(100,32){3} \put(100,16){4}
\put(100,0){5}

\put(118,83){1} \put(134,83){2} \put(150,83){3} \put(166,83){4}
\put(182,83){5}

\put(181,71){\oval(8,8)[br]} \put(165,55){\oval(8,8)[br]}
\put(149,39){\oval(8,8)[br]} \put(133,23){\oval(8,8)[br]}
\put(117,7){\oval(8,8)[br]}

\put(165,71){\oval(8,8)[br]} \put(173,63){\oval(8,8)[tl]}
\put(149,71){\oval(8,8)[br]} \put(157,63){\oval(8,8)[tl]}
\put(149,55){\oval(8,8)[br]} \put(157,47){\oval(8,8)[tl]}
\put(133,39){\oval(8,8)[br]} \put(141,31){\oval(8,8)[tl]}
\put(117,23){\oval(8,8)[br]} \put(125,15){\oval(8,8)[tl]}
\put(117,71){\oval(8,8)[br]} \put(125,63){\oval(8,8)[tl]}


\put (118,51){\line(1,0){6}} \put (121,48){\line(0,1){6}}

\put (118,35){\line(1,0){6}} \put (121,32){\line(0,1){6}}

\put (134,67){\line(1,0){6}} \put (137,64){\line(0,1){6}}

\put (134,51){\line(1,0){6}} \put (137,48){\line(0,1){6}}


\put (108,67){\line(1,0){10}} \put (108,51){\line(1,0){10}} \put
(108,35){\line(1,0){10}} \put (108,19){\line(1,0){10}} \put
(108,3){\line(1,0){10}} \put (124,67){\line(1,0){10}} \put
(124,51){\line(1,0){10}} \put (124,35){\line(1,0){10}} \put
(124,19){\line(1,0){10}} \put (140,67){\line(1,0){10}} \put
(140,51){\line(1,0){10}} \put (140,35){\line(1,0){10}} \put
(156,67){\line(1,0){10}} \put (156,51){\line(1,0){10}} \put
(172,67){\line(1,0){10}}


\put (121,70){\line(0,1){10}} \put (121,54){\line(0,1){10}} \put
(121,38){\line(0,1){10}} \put (121,22){\line(0,1){10}} \put
(121,6){\line(0,1){10}} \put (137,70){\line(0,1){10}} \put
(137,54){\line(0,1){10}} \put (137,38){\line(0,1){10}} \put
(137,22){\line(0,1){10}} \put (153,70){\line(0,1){10}} \put
(153,54){\line(0,1){10}} \put (153,38){\line(0,1){10}} \put
(169,70){\line(0,1){10}} \put (169,54){\line(0,1){10}} \put
(185,70){\line(0,1){10}}



\put(200,64){1} \put(200,48){2} \put(200,32){3} \put(200,16){4}
\put(200,0){5}

\put(218,83){1} \put(234,83){2} \put(250,83){3} \put(266,83){4}
\put(282,83){5}

\put(281,71){\oval(8,8)[br]} \put(265,55){\oval(8,8)[br]}
\put(249,39){\oval(8,8)[br]} \put(233,23){\oval(8,8)[br]}
\put(217,7){\oval(8,8)[br]}

\put(233,71){\oval(8,8)[br]} \put(241,63){\oval(8,8)[tl]}
\put(249,71){\oval(8,8)[br]} \put(257,63){\oval(8,8)[tl]}
\put(249,55){\oval(8,8)[br]} \put(257,47){\oval(8,8)[tl]}
\put(233,39){\oval(8,8)[br]} \put(241,31){\oval(8,8)[tl]}
\put(217,23){\oval(8,8)[br]} \put(225,15){\oval(8,8)[tl]}
\put(217,39){\oval(8,8)[br]} \put(225,31){\oval(8,8)[tl]}


\put (218,67){\line(1,0){6}} \put (221,64){\line(0,1){6}}

\put (218,51){\line(1,0){6}} \put (221,48){\line(0,1){6}}

\put (266,67){\line(1,0){6}} \put (269,64){\line(0,1){6}}

\put (234,51){\line(1,0){6}} \put (237,48){\line(0,1){6}}


\put (208,67){\line(1,0){10}} \put (208,51){\line(1,0){10}} \put
(208,35){\line(1,0){10}} \put (208,19){\line(1,0){10}} \put
(208,3){\line(1,0){10}} \put (224,67){\line(1,0){10}} \put
(224,51){\line(1,0){10}} \put (224,35){\line(1,0){10}} \put
(224,19){\line(1,0){10}} \put (240,67){\line(1,0){10}} \put
(240,51){\line(1,0){10}} \put (240,35){\line(1,0){10}} \put
(256,67){\line(1,0){10}} \put (256,51){\line(1,0){10}} \put
(272,67){\line(1,0){10}}


\put (221,70){\line(0,1){10}} \put (221,54){\line(0,1){10}} \put
(221,38){\line(0,1){10}} \put (221,22){\line(0,1){10}} \put
(221,6){\line(0,1){10}} \put (237,70){\line(0,1){10}} \put
(237,54){\line(0,1){10}} \put (237,38){\line(0,1){10}} \put
(237,22){\line(0,1){10}} \put (253,70){\line(0,1){10}} \put
(253,54){\line(0,1){10}} \put (253,38){\line(0,1){10}} \put
(269,70){\line(0,1){10}} \put (269,54){\line(0,1){10}} \put
(285,70){\line(0,1){10}}
\end{picture}

\caption{ \protect\small Diagrams that are given by
$\{(2,1),(2,2),(2,3),(3,1),(3,2)\}$, $\{(1,2),(2,1),(2,2),(3,1)\}$ and
$\{(1,1),(1,4),(2,1),(2,2)\}$.  }
\label{fig:graphs}
\end{figure}
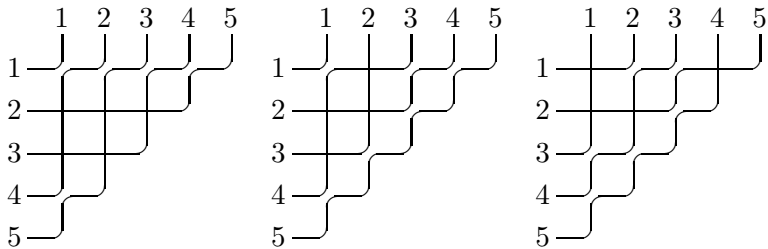

Consider a finite subset~$R$ of $\{1, \ldots,n\}\times \{1\ldots,n\}$ and
think of it as a network of pipes, which intersect at each $(i,j)\in R$
and do not intersect otherwise.  Such subsets are called {\em diagrams}\/
(but are also known as {\em pipe dreams}\/).  For example, see
Figure~\ref{fig:graphs}.

Associate to each diagram $R$ the permutation $w_R \in S_n$ such that the
pipe entering row~$i$ exits column~$w_R(i)$.  For example, the
permutations associated to the diagrams from Figure \ref{fig:graphs} are
$15423$, $14235$ and $21534$.  The diagram~$R$ is an {\em rc-graph}\/ (or
{\em reduced pipe dream}\/) if no two strands intersect twice.  The first
and the third diagrams from Figure~\ref{fig:graphs} are rc-graphs, while
the second one is not.  These diagrams were originally introduced in
\cite{FKyangBax}.  They index the monomials in Schubert polynomials the
same way that semistandard Young tableaux index monomials in Schur
polynomials; see \cite{BJS}, \cite{FSnilCoxeter} and \cite{FKyangBax} for
details.

For an rc-graph $R$, let
\begin{itemize}
\item
$L_R$ be the coordinate subspace of~$\mn$ consisting of all matrices
whose coordinates $z_{ij}$ are zero for every crossing $(i,j)\in R$;

\item
$F_R$ be the {\em rc-face}\/ of the Gel$'$fand--Cetlin polytope given by
setting $\lambda_{i,j}=\lambda_{i+1,j}$ for each $(i,j)\in R$; and

\item
$T_{\hspace{-.1ex}R}$ be the {\em rc-toric subvariety}\/ of the
Gel$'$fand--Cetlin toric variety with face~$F_R$.
\end{itemize}

Next let us review some geometric ingredients for our main theorem and
its proof.  For a permutation $w \in S_n$, the {\em Schubert
determinantal ideal}\/ $I_w\subseteq \CC[\zz]$, defined by Fulton in
\cite{Ful92}, is generated by all minors of size $1+w_{qp}$ in the top
left $q\times p$ submatrix~$Z_{p\times q}$ of $Z=(z_{ij})$ for all
$q,p$, where $w_{qp}$ is the number of $i\leq q$ such that $w(i)\leq p$.
The {\em matrix Schubert variety}\/ for~$w$ \cite{grobGeom} is by
definition the zero set~$\ol X_w$ of~$I_w$.  We denote by
$X_w\subset\FL$ the {\em Schubert variety} obtained by projecting $\ol
X_w \cap \gln$ to~$\FL$.  This Schubert variety is the closure in~$\FL$
of the $B^+$ orbit through the coset~$B \ol w$, where the permutation
matrix~$\ol w$ has its nonzero entries at~$(i,w(i))$, and $B^+$ is the
Borel group of upper triangular matrices acting on the right of $\FL =
\bgln$.

In general, a flat family degenerating a variety~$Y$ does not induce
degenerations on its subvarieties.  However, Gr\"obner and sagbi
degenerations of~$Y$ are canonically isomorphic to trivial families
over~$\CC^*$.  Any subvariety of~$Y$ defines an (isomorphically trivial)
subfamily over~$\CC^*$, and hence a flat subfamily over all of~$\CC$ by
taking the closure of this subfamily.

\begin{thm} \label{t:main}
The quotient family $\FF = B \GIT (\mn \times \CC)$ induces flat
degenerations of Schubert subvarieties $X_w$ of the complete flag
manifold $\FL = \FF(1)$ to reduced unions
$\bigcup_{w(R)=w} T_{\hspace{-.1ex}R}$ of toric subvarieties of the
Gel\/$'\!$fand--Cetlin toric~variety $\FF(0)$.
\end{thm}
\begin{proof}
It was shown in \cite{grobGeom} that the minors generating~$I_w$
constitute a Gr\"obner basis for any antidiagonal term order, and
hence define a Gr\"obner degeneration of~$\ol X_w$.  Moreover, it
was shown that any such degeneration---in particular (by
Lemma~\ref{l:weights}) the one given by~$\omega$---degenerates the
matrix Schubert variety $\ol X_w$ to the reduced union
$\bigcup_{w(R)=w}L_R$ of rc-subspaces for~$w$ inside~$\mn$.  The
image of the total space of this Gr\"obner degeneration under the
family of Pl\"ucker embeddings given by coordinates~(\ref{q})
equals our family~$\FF$ by Theorem~\ref{t:degeneration}.  On the
other hand, the closure of the image of an rc-subspace $L_R$ under
the degenerated Pl\"ucker map obtained by setting $t=0$
in~(\ref{q}) equals the corresponding rc-toric
subvariety~$T_{\hspace{-.2ex}R}$ by definition.%
\end{proof}

The argument in the proof can be summarized as: the GIT quotient by~$B$
of the Gr\"obner degeneration in \cite{grobGeom} equals the sagbi
degeneration in Theorem~\ref{t:main}.

\begin{remark}
The Schubert variety $X_w \subseteq \FL$ equals the intersection of the
embedded subvariety $\FL \subset \prod\PP(\bigwedge^{\!k}\CC^n)$ with a
set of hyperplanes, one hyperplane $p_I=0$ for each subset
$I=\{i_1,\ldots,i_k\}$ satisfying $k>w_{k,i_k}$, where $w_{k,i_k}$ is the
number of $i\leq k$ with $w(i)\leq i_k$.  Intersecting the family~$\FF$
with the same set of hyperplanes produces the degeneration of~$X_w$ in
Theorem~\ref{t:main}.%
\end{remark}

\begin{remark}
Using the involution on~$\FL$ that switches the Pl\"ucker coordinate
$p_I$ with $p_{\bar I}$, where $\bar I$ is the complement of~$I$, it can
be shown that opposite Schubert varieties degenerate to unions of toric
subvarieties associated to opposite rc-walls of the Gel$'$fand--Cetlin
polytope.%
\end{remark}

\section{Gel\cprime fand--Cetlin decomposition}

\label{sec:decomp}

This final section gives a geometric construction of the
Gel$'$fand--Cetlin basis of an irreducible $\gln$ representation, by
extending the \mbox{Borel--Weil} construction to the whole family~$\FF$.
Our proof logically depends only on the Borel--Weil theorem.  To
introduce notation, we begin by reviewing the construction of Gel$'$fand
and Cetlin~\cite{GC50}.

For a dominant weight $\lambda$ of~$\gln$, which by definition is a
decreasing sequence $(\lambda_1>\cdots>\lambda_n)$ of positive integers,
let $V^\lambda$ be the irreducible representation of $\gln$ with highest
weight~$\lambda$.  For $n\geq i\geq 1$ identify $\gli$ with the subgroup
of $\gln$ sitting in the bottom right \mbox{$i\times i$} corner.  As a
$\glnone$ representation, $V^\lambda$ breaks up into a direct sum of
irreducible components
\begin{eqnarray} \label{eq:weyl}
  V^\lambda &=& \bigoplus_{\mu\prec\lambda} V^\mu,
\end{eqnarray}
where the dominant weight $\mu=(\mu_1>\cdots>\mu_{n-1})$ of $\glnone$
satisfies $\mu\prec\lambda$ if
$$
  \lambda_1\geq\mu_1\geq\lambda_2\geq\mu_2\geq\lambda_3\geq\cdots\geq
  \mu_{n-1}\geq\lambda_n,
$$
so $\mu$ interpolates between~$\lambda$.  Iterating defines
\emph{partial Gel$'$fand--Cetlin decompositions}
\begin{eqnarray} \label{eq:partialGC}
  V^\lambda &=& \bigoplus_{\Lambda_{i}} V^{\Lambda_{i}}
\end{eqnarray}
of $V^\lambda$ into irreducible components for the action
of~$\gli$, where $\Lambda_{i}$ runs through chains
$(\lambda\succ\lambda^{n-1}\succ\cdots\succ\lambda^i)$ with
$\lambda^j$ being a weight of~$G_{\!}L_j$.  Gel$'$fand and Cetlin
studied the decomposition of $V^\lambda$ as a direct sum of
one-dimensional subspaces~$V^\Lambda$, one for each chain
$\Lambda=(\lambda \succ \lambda^{n-1}\succ\cdots\succ \lambda^1)$.
By definition, $\Lambda$ lies in~$\Pi_\lambda$, the set of integer
Gel$'$fand--Cetlin patterns for~$\lambda$.  Hence we have the {\em
Gel\/$'\!$fand--Cetlin decomposition}
\begin{eqnarray} \label{GC}
  V^\lambda &=& \bigoplus_{\Lambda\in \Pi_\lambda} V^\Lambda.
\end{eqnarray}

For a weight $\lambda$, consider the very ample line bundle
$\LL=\OO_\FF (\aa)$ from Section~\ref{sec:GC} over the
family~$\FF$. Let $\LL^\lambda_\tau$ be the restriction of this
line bundle to the fiber over $\tau\in \CC$.  The Borel--Weil
theorem states that the representation $V^\lambda$ with highest
weight~$\lambda$ is isomorphic to the space of algebraic sections
of~$\LL_1^\lambda$ as a representation of~$\gln$:
\begin{eqnarray} \label{BW}
  V^\lambda &=& \Gamma(\LL_1^\lambda).
\end{eqnarray}
At the same time, we have already seen that the space
$\Gamma(\LL^\lambda_0)$ of sections over the toric
variety~$\FF(0)$ carries an action of a torus $\TT$ of dimension
$\binom{n}{2}$ under which $\Gamma(\LL^\lambda_0)$ decomposes into
one-dimensional weight spaces.  Think of~$\TT$ as the product of
one-dimensional tori~$T'_{ij}$ for $i+j \leq n$, each of which
acts on~$\mn$ by scaling the $(i,j)$ entry, which lies strictly
above the main antidiagonal.  The weight spaces for the action
of~$\TT$ on $\Gamma(\LL^\lambda_0)$ are indexed by the set~$\UL$
from~(\ref{UL}) in Section~\ref{sec:GC}.  In other words,
\begin{eqnarray} \label{LL1}
  \Gamma(\LL^\lambda_0) &=& \bigoplus_{A \in \UL}\CC^A,
\end{eqnarray}
where $\CC^A$ is a complex line on which $\TT$ acts with weight~$A$.

Let $T_{ij}$ be the one-dimensional torus scaling simultaneously the
entries of an $n\times n$ matrix in row~$i$, between column~\mbox{$j+1$}
and the antidiagonal column~\mbox{$n+1-i$}.  Then $\TT$ can be thought
of as the product of all tori $T_{ij}$ with $i+j\leq n$.  Under this
direct product decomposition of~$\TT$, the discussion in the proof of
Proposition~\ref{prop:GC-polytope} identifying $\UL$ with~$\Pi_\lambda$
implies the weight space decomposition
\begin{eqnarray} \label{LL}
  \Gamma(\LL^\lambda_0) &=& \bigoplus_{\Lambda\in
  \Pi_\lambda}\CC^\Lambda,
\end{eqnarray}
The weight space decompositions (\ref{LL1}) and~(\ref{LL}) are of course
the same, and the two indexings correspond to two different choices of
bases of the weight lattice of $\TT$.

The family~$\FF$ is projective over the affine complex line
$\spec(\CC[t])$, so the algebraic sections $\Gamma(\LL^\lambda)$ form a
finitely generated module over the coordinate ring $\CC[t]$ of the base.
The localization $\Gamma(\LL^\lambda)\otimes_{\CC[t]} \CC[t^{\pm1}]$ is
a finitely generated free module over the coordinate ring
$\CC[t^{\pm1}]$ of the complement of $0 \in \CC$, by triviality of the
family~$\FF$ outside the fiber over~$0$, and invariance of the vector
space dimension of~$\Gamma(\LL^\lambda_\tau)$ as a function of~$\tau$.
On the other hand, $\dim_\CC \Gamma(\LL^\lambda_0) =
\dim_\CC(\LL^\lambda_\tau)$ for $\tau \neq 0$; indeed, both dimensions
equal the number of lattice points in the Gel$'$fand--Cetlin
polytope~$P_\lambda$, by~(\ref{GC}), (\ref{BW}), and~(\ref{LL}).  Hence
we get the following.

\begin{prop} \label{p:Phi}
$\Gamma(\LL^\lambda)$ is free over~$\CC[t]$, and possesses a
$\CC[t]$-basis of sections each of which is equivariant for the action
of the $n$-dimensional diagonal torus in~$\gln$.  Restricting this basis
to the\/~$0$ and\/~$1$ fibers results in a torus-equivariant isomorphism
$$
\begin{array}{rccl}
  \Phi:& \Gamma(\LL^\lambda_1) = V^\lambda &\too& \displaystyle
  \bigoplus_{\Lambda\in \Pi_\lambda}\CC^\Lambda = \Gamma(\LL^\lambda_0).
\end{array}
$$
\end{prop}


To understand the map $\Phi$ in terms of the inductive construction of
Gel$'$fand and Cetlin, we present a construction of an $n$-parameter
family extending~$\FF$.
Let $\omega_i$ be the $n\times n$ matrix whose $i^\th$~column has
entries $\omega_{1i},\ldots,\omega_{ni}$ and whose other columns
are zero (the integers~$\omega_{ij}$ are defined
in~(\ref{eq:matrix})). Write $T(\omega_i)$ for the one parameter
subgroup of the torus of~$B^n$ associated to~$\omega_i$, and
denote by $\tilde\tau_i$ the element of $T(\omega_i)$
corresponding to the complex number~$\tau_i$.  Define the family
\begin{eqnarray*}
  B(\tau_1,\cdots,\tau_n) &=& \tilde\tau_1^{-1} \cdots \tilde\tau_n^{-1}
  \BD \tilde\tau_1\cdots \tilde\tau_n
\end{eqnarray*}
of subgroups of~$B^n$.  This family extends to zero values of $\tau_i$
and defines an action of~$B$ on $\mn\times \CC^n$ as in
Lemma~\ref{l:tau}.
\begin{defn}
The $n$-parameter family $\tilde \FF = B \GIT (\mn\times \CC^n)$ is the
{\em degeneration in stages}.  Denote by~$\tilde\FF_i$ its fiber over
the point $(0,\ldots,0,1\ldots,1)$ with $i$ entries equal~to~$1$.
\end{defn}

Observe that $\tilde\FF_n = \FL$ is the flag manifold, and~$\tilde\FF_0$
is the Gel$'$fand--Cetlin toric variety by Theorem~\ref{t:main}.

Let $\TT_i$ be the torus \mbox{$T_{n-1} \times \cdots \times T_i$} with
$T_j = T_{1,n-j} \times \cdots \times T_{j,n-j}$, and set~\mbox{$\GG_i =
G_{\!}L_{i}$}, thought of in the bottom right corner again.  For the
$n=5$ example of~$\TT_i$, each torus $T_j$ scales the entries by its
one-parameter subgroups $T_{ij}$ in the \mbox{indicated locations}:
$$
\def\*#1#2{\makebox[2.5ex][l]{$T_{#1#2}$}}
\def\ct#1#2{\multicolumn{1}{|@{\,}l@{}|}{\makebox[2ex][l]{$\*#1#2$}}}
\def\lt#1#2{\multicolumn{1}{|@{\,}l@{}}{\makebox[2ex][l]{$\*#1#2$}}}
\def\ro{\multicolumn{1}{c|}{}}
T_4\colon\ \tinyrc{\begin{array}{|@{}lcccc|}\hline
    &\lt11 &   &   &\ro\\\cline{2-5}
    & \lt21&   &\ro&   \\\cline{2-4}
    &\lt31 &\ro&   &   \\\cline{2-3}
    & \ct41&   &   &   \\ \cline{2-2}
\ \ &\ \   &   &   &\, \\\hline
\end{array}}
\:,\quad
T_3\colon\ \tinyrc{\begin{array}{|c@{}lccc|}\hline
    &    &\lt12&   &\ro\\\cline{3-5}
    &    &\lt22&\ro&   \\\cline{3-4}
    &    &\ct32&   &   \\\cline{3-3}
    &    &     &   &   \\
\ \ &\ \ &\ \  &   &\, \\\hline
\end{array}}
\:,\quad
T_2\colon\ \tinyrc{\begin{array}{|c@{}lccc|}\hline
    &  &  &\lt13&\ro\\\cline{4-5}
    &  &  &\ct23&   \\\cline{4-4}
    &  &  &     &   \\
    &  &  &     &   \\
\ \ &\ &\ &\ \, &\, \\\hline
\end{array}}
\:,\quad
T_1\colon\ \tinyrc{\begin{array}{|c@{}lccc|}\hline
    &    &    &  &\ct14\\\cline{5-5}
    &    &    &  &     \\
    &    &    &  &     \\
    &    &    &  &     \\
\ \,&\ \,&\,  &\ &\ \: \\\hline
\end{array}}
$$
Each fiber $\tilde\FF_{i}$ carries the action of the group $\GG_i \times
\TT_i$.  Moreover, for each~$i$, a~subfamily of $\tilde \FF$ degenerates
$\tilde\FF_{i}$ to $\tilde\FF_{i-1}$ flatly and $\GG_{i-1}\times \TT_i$
invariantly.

Associate to each dominant weight~$\lambda$ the very ample line bundle
$\tilde \LL^\lambda$ over $\tilde \FF$, as we did for the
family~$\FF$. Then, as in~(\ref{PL}) for~$\FF$, we get an embedding
of~$\tilde\FF$ into $\PP_\lambda \times \CC^n$.  Let $V^\lambda_i$ be
the space of algebraic sections of the restriction of the line
bundle~$\tilde\LL^\lambda$ to~$\tilde\FF_{i}$, treated as a
representation of $\GG_i \times \TT_i$.  Since $\tilde\FF_{i}$
degenerates to $\tilde\FF_{i-1}$ flatly and $\GG_{i-1} \times \TT_i$
invariantly, there are~$\GG_{i-1} \times \TT_i$ invariant isomorphisms
$\Phi_{i}: V^\lambda_i \to V^\lambda_{i-1}$ for $i = 1,\ldots,n$. These
are analogous to the isomorphism~$\Phi$ from Proposition~\ref{p:Phi},
and constructed geometrically in the same way. In fact $\Phi$ equals the
composition $\Phi_1 \circ \cdots \circ \Phi_n$, or equivalently $\Phi:
V^\lambda=V^\lambda_n \stackrel{\Phi_n}\too V^\lambda_{n-1}
\stackrel{\Phi_{n-1}}\too \cdots \stackrel{\Phi_1}\too V^\lambda_0$.

In what follows, we write $\Pi_\lambda(i)$ for the set of integer
patterns
\begin{eqnarray*}
  \Lambda_i &=&
  (\lambda\succ\lambda^{n-1}\succ\cdots\succ\lambda^{i}),
\end{eqnarray*}
with $\lambda^j$ being a weakly decreasing sequence of\/~$j$ nonnegative
integers.  For each pattern $\Lambda_i \in \Pi_\lambda(i)$, let
$V_i^{\Lambda_i}$ be the irreducible representation of\/~$\GG_i$ with
highest weight~$\lambda^i$, and declare the torus $\TT_i$ to act on
every vector in $V_i^{\Lambda_i}$ with weight~$\Lambda_i$.

\begin{thm} \label{thm:basis}
The sections $V^\lambda_i$ of the line bundle~$\tilde\LL^\lambda$
over the fiber $\tilde\FF_{i}$ of the degeneration
in stages decomposes into irreducible components for\/ \mbox{$\GG_i
\times \TT_i$} as
\begin{eqnarray} \label{Vi}
  V_i^\lambda &=& \bigoplus_{\Lambda_i \in \Pi_\lambda(i)}
  V_i^{\Lambda_i},
\\\nonumber
  \hbox{so\quad} V^\lambda &=& \bigoplus_{\Lambda_i \in \Pi_\lambda(i)}
  \Phi_n^{-1} \circ \cdots \circ \Phi_{i+1}^{-1}(V_i^{\Lambda_i})
\end{eqnarray}
is a partial Gel\/$'\!$fand--Cetlin decomposition.  Thus
the Gel\/$'\!$fand--Cetlin decomposition~is
\begin{eqnarray*}
  V^\lambda &=& \bigoplus_{\Lambda \in \Pi_\lambda}
  \Phi^{-1}(\CC^\Lambda).
\end{eqnarray*}
\end{thm}
\begin{proof}
It is enough to assume~(\ref{Vi}) is proved for~$i$, and then prove it
for~\mbox{$i-1$}.  Let $\Lambda_{i-1} \mapsto \Lambda_i$ be the map
$\Pi_\lambda(i-1) \to \Pi_\lambda(i)$ forgetting $\lambda^{i-1}$. Since
$\Phi_{i}$ is $\GG_{i-1}$ equivariant, we get a decomposition of
$V_{i-1}^\lambda$ into irreducible components under~$\GG_{i-1}$:
\begin{eqnarray*}
  V_{i-1}^\lambda &=& \bigoplus_{\!\Lambda_{i-1} \in \Pi_\lambda(i-1)\,}
  \bigoplus_{\,\Lambda_{i-1} \mapsto \Lambda_i} V^{\lambda^{i-1}\!}.
\end{eqnarray*}
It remains to show that
$T_{i-1}$ acts on each irreducible $V^{\lambda^{i-1}}$ with
weight~$\lambda^{i-1}$.

After the identification of $V^\lambda$ with $\Gamma(\LL^\lambda_1)$,
every highest weight vector of the $\GG_{i-1}$ action on $V^\lambda$ can
be thought of as a monomial $\prod p_I$ in Pl\"ucker coordinates for
subsets~$I$ whose columns in the range $n-i+2,\ldots,n$ are left
justified.  (These are the monomials invariant with respect to the right
action of $U^+_{i-1}$, the upper triangular matrices inside $\GG_{i-1}$
with $1$'s on the diagonal.)  Now simply note that the weight
of~$T_{i-1}$ on such a monomial coincides with the weight of the
diagonal torus in~$\GG_{i-1}$.%
\end{proof}

\subsection*{Acknowledgements}
Both authors thank Allen Knutson for numerous useful discussions.  In
particular, it was he who first suggested to us that the flag manifold
could be degenerated to the Gel$'$fand--Cetlin toric variety, in stages
as well.  Thanks also to Peter Littelmann, for pointing out helpful
references.


\end{document}